\documentclass[11pt,draft]{article}
\usepackage[cp1250]{inputenc}
\usepackage{latexsym,amssymb,amsmath}

\title{On Frattini subloops and normalizers of commutative Moufang loops}
\date{23.04.2008}
\author{N. I. Sandu}
\begin{document}
\maketitle

\begin{abstract}
Let $L$ be a commutative Moufang loop (CML) with multiplication
group $\frak M$, and let $\frak F(L)$, $\frak F(\frak M)$ be the
Frattini subgroup and Frattini subgroup of $L$ and $\frak M$
respectively. It is proved that $\frak F(L) = L$ if and only if
$\frak F(\frak M) = \frak M$ and is described the structure of
this CLM. Constructively it is defined the notion of normalizer
for subloops in CML. Using this it is proved that if $\frak F(L)
\neq L$ then $L$ satisfies the normalizer condition and that any
divisible subgroup of $\frak M$ is an abelian group and serves as
a direct factor for $\frak M$.

\ Classification: {\ 20N05}

\end{abstract}

It is known that in many classes of algebras the Frattini
subalgebras essential define the structure of these algebras. In
this paper this dependance is considered in the class of
commutative Moufang loops (CML) and its multiplication groups. Let
$L$ be a CML with multiplication group $\frak M$, let $\frak F(L)$
and $\frak F(\frak M)$ denote the Frattini subloop and Frattini
subgroup of $L$ and $\frak M$ respectively. It is proved that
$\frak F(L) = L$ if and only if $\frak F(\frak M) = \frak M$ and
the structure of this CML and groups is described. In particular,
if $L$ has the exponent $3$ then $\frak F(L) = L$ if and only if
$L = L^{\prime}$, where $L^{\prime}$ denotes the associator
subloop of $L$ (Theorem 1). The existence of CML with $L^{\prime}
= L$ is proved in [1].

Constructively the normalizer $N_L(H)$  is defined  for subloop
$H$ of commu\-tative Moufang loop $L$ which, in general, has the
same destina\-tion as a normalizer for subgroups. The normalizer
$N_L(H)$ is the unique maximal subloop of $L$ such that $H$ is
normal in $N_L(H)$. By analogy with the group theory it is defined
the notion of CML with normalizer condition. Thus is then, when
every proper subloop of CML differs from its normalizer.
Essentially using the Theorem 1, it is proved that if a CML $L$
satisfy the inequality $\frak F(L) \neq L$ then $L$ satisfies the
normalizer condition. It is proved also that for multiplication
groups of CML an analogous situation   is false. There exists a
CML $L$ with multiplication group $\frak M$ such that $\frak
F(\frak M) \neq \frak M$. but $\frak M$ does not satisfy the
normalizer condition.

Again, essentially  using the Theorem 1 it is proved that every
divisible subgroup of multiplication group $\frak M$ of all CML is
an abelian group and serves as a direct factor for $\frak M$
(Theorem 2). A similar  result for divisible subloops of CML is
proved in [2]. We denote that in general case the Theorem 2 is
false. In [3, Theorem 2.7 and Example 2.2] there is an example of
divisible non-periodic and non-abelian $ZA$-group.

\section{Preliminaries}

Let us bring some notions and results on the loop theory  from [4,
5].

The \textit{multiplication group} $\frak M(L)$ of a loop $L$ is
the group generated by all  \textit{translations} $L(x)$, $R(x)$,
where $L(x)y = xy$, $R(x)y = yx$. The subgroup $\frak I(L)$ of
group $\frak M(L)$, generated by all  \textit{inner mappings}
$L(x,y) = L^{-1}(xy)L(x)L(y)$, $R(x,y) = R^{-1}(xy)R(y)R(x)$,
$T(x) = L^{-1}(x)R(x)$  is cal\-led the \textit{inner map\-ping
group} of  loop $L$. The subloop $H$ of  loop $L$ is called
\textit{normal} in $L$, if $\frak I(L)H = H$. The set of all
elements $x \in L$ which commute and associate with all elements
of $L$, so that for all $a, b \in L$ $ax = xa, ab\cdot x = a\cdot
bx, ax\cdot b = a\cdot xb, xa\cdot b = x\cdot ab$ is a normal
subloop $Z(L)$ of $L$, its \textit{centre}.
\smallskip\\

\textbf{Lemma 1} [4, page 62]. \textit{Let $H$ be a normal subloop
of loop $L$ with multiplication group $\frak M$. Then $\frak
M(L/H) \cong \frak M/H^{\ast}$ where $H^{\ast} = \break \{\alpha
\in \frak M \vert (\alpha x)H = xH \quad \forall x \in L\}$.
Conversely, every normal subgroup $\frak N$ of $\frak M$
determines a normal subloop $H = \frak N1 = \{\alpha1 \vert \alpha
\in \frak N\}$ of $L$ and $\frak N \subseteq H^{\ast}$.}
\smallskip\\

\textbf{Proposition 1.} \textit{Let $(L, \cdot, 1)$ be a loop with
centre $Z(L)$, let $\frak M$ be its multiplication group with
centre $Z(\frak M)$ and let $\tilde{Z}(L) = \{\varphi 1 \vert
\varphi \in Z(\frak M)\}$, $\tilde{Z}(\frak M) = \{L(\varphi 1)
\vert \varphi \in Z(\frak M)\}$, $\overline{Z}(\frak M) = \{L(a)
\vert a \in Z(L)\}$. Then $\overline{Z}(L) = Z(L) \cong
\overline{Z}(\frak M) = \tilde{Z}(\frak M) = Z(\frak M).$}
\smallskip\\

\textbf{Proof.} Let $a \in Z(Q)$ and $x, y \in L$. Then $R(a) =
L(a)$, $a\cdot a^{-1}x = x$, $L(a)L(a^{-1})x = x$ $L(a^{-1}) =
L^{-1}(a)$ and $a\cdot xy = ax\cdot y$, $L(a)L(y)x = L(y)L(a)x$,
$L(a)L(y) = L(y)L(a)$. Similarly, for $a^{-1} \in Z(L)$ we
obtained that $L(a^{-1})R(y) = R(y)L(a^{-1})$. Then
 $(L(a^{-1})R(y))^{-1} = (R(y)L(a^{-1}))^{-1}$,
$R^{-1}(y)L^{-1}(a^{-1}) = L^{-1}(a^{-1})R^{-1}(y)$,
$R^{-1}(y)L(a) = L(a)R^{-1}(y)$. Analogo\-usly, from $yx\cdot a =
y\cdot xa$ and $R(a) = L(a)$ we get $L(a)L(y) = L(y)L(a)$,
$L(a)L^{-1}(y) = L(y)^{-1}L(a)$. Then from definition of group
$\frak M(L)$ it follows that $L(a) \in Z(\frak M)$. Similarly, for
$a^{-1} \in Z(L)$ we get that $L(a^{-1}) = L^{-1}(a) \in Z(\frak
M)$. We also  have  $L(a)L(b) = L(ab)$ for $a, b \in Z(L)$. Then
the set $\overline{Z}(\frak M)$ is a subgroup of $\frak M$,
$\overline{Z}(\frak M)\subseteq Z(\frak M)$ and the isomorphism
$Z(L) \cong \overline{Z}(\frak M)$, defined by $u \rightarrow
L(u)$, $u^{-1} \rightarrow L^{-1}(u)$, $u \in L$, follows from the
equality $L(a)L(b) = L(ab)$.

If $\varphi \in Z(\frak M)$, then  $\varphi L(x) = L(x)\varphi$,
$\varphi L(x)y = L(x)\varphi y,  \varphi (xy) = x\cdot \varphi y$
and $\varphi R(x) = R(x)\varphi$, $\varphi R(x)y = R(x)\varphi y$,
$\varphi(yx) = \varphi y\cdot x$ for any $x, y \in L$. Hence
$\varphi (xy) = x\cdot \varphi y$ and $\varphi(yx) = \varphi
y\cdot x$. Let $y = 1$. Then $\varphi x = x\cdot \varphi 1$,
$\varphi x = \varphi 1 \cdot x$, i.e. $x\cdot \varphi 1 = \varphi
1\cdot x$.  Now, using the equality $\varphi(xy) = x\cdot \varphi
y$ we obtain that $xy\cdot \varphi 1 = \varphi (xy \cdot 1) =
\varphi (xy) = x\cdot\varphi y = x\cdot\varphi (y\cdot 1) =
x(y\cdot \varphi 1)$ and using the equality $\varphi(yx) = \varphi
y\cdot x$ we obtain that $\varphi 1 \cdot xy = \varphi(1 \cdot xy)
= \varphi(xy) = \varphi x \cdot y = \varphi(1\cdot x)y = (\varphi
1 \cdot x)y$. Hence, if $\varphi \in Z(\frak M)$ then $\varphi 1
\in Z(L)$, i.e. $\tilde{Z}(L) \subseteq Z(L)$. Conversely, let $a
\in Z(L)$. Then $L(a) \in Z(\frak M)$ and $a = L(a)1 \in
\tilde{Z}(L)$. Hence $Z(L) \subseteq \tilde{Z}(L)$ Consequently,
$Z(L) = \tilde{Z}(L)$ and therefore $\overline{Z}(\frak M) =
\tilde{Z}(\frak M)$.

Let $\frak I(L)$ be the inner mapping group of $\frak M$. In the
proof of Lemma IV.1.2 from [4] it is shown that each element
$\alpha \in \frak M$ has the form $\alpha = L(\alpha 1)\theta$
where $\theta \in \frak I(L)$; moreover $\alpha \in \frak I(L)$ if
and only if $L(\alpha 1) = e$ where $e$ is the unit of $\frak M$.
Let $\frak J(L) = \frak I(L) \cap Z(\frak M)$. If $\alpha \in
Z(\frak M)$ then, by the cases considered above, $L(\alpha 1) \in
\overline{Z}(\frak M)$. Then $\theta \in \frak{J}(L)$. The
subgroups $\overline{Z}(\frak M) \subseteq Z(\frak M)$ and $\frak
J(L) \subseteq Z(\frak M)$ are normal in $\frak M$. As
$\overline{Z}(\frak M) \cap \frak{J}(L) = \varepsilon$ then
$Z(\frak M) = \overline{Z}(\frak M) \times \frak{J}(L)$. By Lemma
1 $\frak J(L)1 = 1$ is a normal subloop of $L$ and $\frak J(L)
\subseteq 1^{\ast}$ where $1^{\ast} = \{\alpha \in \frak M \vert
(\alpha x)1 = x1 \quad \forall x \in L\}$. But $1^{\ast} = e$,
hence $\frak J(L) = e$ and $Z(\frak M) = \overline{Z}(\frak M)$.
This completes the proof of Proposition 1.
\smallskip\\

A system $\Sigma$ of subloops of loop $L$ will be called
\textit{subnormal system }, if: 1) it contains $1$ and $L$; 2) is
linearly ordered by inclusion, i. e. for all $A, B$ from $\Sigma$
either $A \subseteq B$, or $B \subseteq A$; 3) is closed regarding
the unions and intersections, in particular, together with each $A
\neq L$ contains the intersection $A^{\sharp}$ of all $H \in
\Sigma$ with the condition $H \supset A$ and together with each $B
\neq 1$ contains the union $B^{\flat}$ of all $H \in \Sigma$ with
the condition $H \subset B$; 4) satisfies the condition: $A$ is
normal in $A^{\natural}$ for all $A \in \Sigma$, $A \neq L$.

A system $\Sigma$ is called \textit{ascending} (respect.
\textit{descending}) if $A^{\sharp} \neq A$ (respect. $B^{\flat}
\neq B$) for all $A \in \Sigma$, $A \neq L$ (respect. $B \in
\Sigma$, $B \neq 1$) and is called \textit{normal} if the subloops
$A \in \Sigma$ are normal in $L$.

A loop $L$ may be called an \textit{$SD$-loop} if it has a
descending subnormal system $\Sigma$, such that the quotient loops
$A^{\sharp}/A$ are abelian groups for all $A \in \Sigma$, $A \neq
L$. If a loop $L$ has an ascending normal system, such that
$A^{\sharp}/A \subseteq Z(L/A)$ for  all $A \in \Sigma$, $A \neq
L$ then $L$ is called a \textit{$ZA$-loop}.

If the upper central series of the $ZA$-loop have a finite length,
then the loop is called \textit{centrally nilpotent}. The least of
such a length is called the \textit{class} of the central
nilpotentcy. If the loop $L$ is centrally nilpotent of class $k$
then the  upper central series of $L$ has the form

$$1 = Z_0(L) \subset Z_1(L) \subset \cdots \subset Z_k(L) = L,\eqno{(1)}$$
where $Z_1(L) = Z(L)$, $Z_{i+1}(L)/Z_i(L) = Z(L/Z_i(L))$.

A \textit{commutative Moufang loop} (\textit{CML}) is
characterized by the identity $ x^2\cdot yz = xy\cdot xz.$ The
\textit{associator} $(a,b,c)$ of the elements $a, b, c$ of the CML
$Q$ are defined by the equality $ab\cdot c = (a\cdot bc)(a,b,c)$.
The identities

$$L(x,y)z = z(z,y,x), \eqno{(2)}$$

$$(x,y,z) = (y,z,x) = (y^{-1},x,z) = (y,x,z)^{-1}, \eqno{(3)}$$

$$(xy,u,v) = (x,u,v)((x,u,v),x,y)(y,u,v)((y,u,v),y,x)
\eqno{(4)}$$ hold  in the CML.

 If $A, B, C$ are subsets of CML $L$,
$(A,B,C)$ denotes the set of all associ\-ators $(a,b,c)$, $a \in
A$, $b \in B$, $c \in C$. If $A = B = C = L$, then the normal
subloop $L^{\prime} = (L,L,L)$ is called the \textit{associator
subloop} of CML $L$.
\smallskip\\

\textbf{Lemma 2} [4]. \textit{Let $L$ be a CML with centre $Z(L)$
and let $a \in L$. Then $a^3\in Z(L)$.}
\smallskip\\

\textbf{Lemma 3.} \textit{If $Z_2(L) \neq Z_1(L)$ for a CML $L$
then $L^{\prime} \neq L$.}
\smallskip\\

\textbf{Proof.} If $z \in Z_2(L) \backslash Z_1(L)$ then
$((x,y,z),u,v) = 1$ for all $x, y, u, v \in L$ and there exist
elements $x_0, y_0 \in L$ such that $(x_0,y_0,z) \neq 1$. From (4)
it follows that $(uv,y_0,z) = (u,y_0,z)(v,y_0,z)$, which show
that the mapping $\varphi: u \rightarrow (u,y_0,z)$ is a
homomorphism of $L$ into $Z_1(L)$. The centre $Z_1(L)$ is an
associative subloop and as $(x_0,y_0,z) \neq 1$ then $L^{\prime}
\subseteq \text{ker}\varphi$ and $L/\text{ker}\varphi$ is
non-unitary. Hence $L^{\prime} \neq L$, as required.
\smallskip\\

\textbf{Lemma 4.} \textit{Let $L$ be a CML with multiplication
group $\frak M$, let $L^{\prime}$ be the associator subloop of $L$
and let $\frak M^{\prime}$ be the commutator subgroup of $\frak
M$. Then $L^{\prime} \subseteq \frak F(L)$ and $\frak M^{\prime}
\subseteq \frak F(\frak M)$.}
\smallskip\\

\textbf{Proof}. The inclusion $L^{\prime} \subseteq \frak F(L)$ is
proved in [6]. As the group $\frak M$ is locally nilpotent [2],
then the proof of inclusion $\frak M^{\prime} \subseteq \frak
F(\frak M)$ can be found, for example, in [7].

\section{Frattini subloops}

If $S, T, \ldots,$ are subsets of elements of a loop $L$, let $<S,
T, \ldots>$ denote the subloop of $L$ generated by $S, T, \ldots.$
An element $x$ of a loop $L$ is a \textit{non-generator} of $L$
if, for every sunset $S$ of $L$, $<x, S> = L$ implies $<S> = L$.
The non-generators of $L$ form the \textit{Frattini subloop},
$\frak F(L)$, of $L$. If $L$ has at least one maximal proper
subloop, then $\frak F(L)$ is the intersection of all maximal
proper subloops of $L$. In the contrary case, $\frak F(L) = L$
[4].
\smallskip\\

\textbf{Lemma 5}. \textit{Let $\theta$ be a homomorphism of the
loop $L$ upon a loop and let $\frak F(L) = L$. Then $\frak
F(\theta L) = \theta L$.}
\smallskip\\

\textbf{Proof.} In [4] it is proved that if $\varphi$ is a
homomorphism of the loop $L$ upon a loop, then $\varphi(\frak
F(L)) \subseteq \frak F(\varphi(L))$. In our case we have $\theta
L = \theta(\frak F(L)) \subseteq \frak F(\theta L) \subseteq
\theta L$. Hence $\frak F(\theta L)) = \theta L$, as required.
\smallskip\\

\textbf{Lemma 6.} \textit{For a CML $L$ with multiplication group
$\frak M$ the following statements are equivalent: 1) $\frak
F(\frak M) = \frak M$; 2) $\frak F(L) = L$.}
\smallskip\\

\textbf{Proof.} $1) \Rightarrow 2)$. Let $\frak F(\frak M) = \frak
M$ and we assume that $\frak F(L) \neq L$. Then $L$ has at least
one maximal proper subloop $H$ and $\frak F(L)$ is the
intersection of all such subloops. By Lemma 4 the associator
subloop $L^{\prime}$ lies in $\frak F(L)$. Hence $H$ is a normal
subloop of $L$ and quotient loop $L/H$ is a cyclic group of prime
order $p$. Then $\frak M(L/H)$ is a cyclic group of order $p$ also
and by Lemma 1 $\frak M/H^{\ast}$ is a cyclic group of order $p$.
Consequently, $H^{\ast}$ is a maximal proper subgroup of $\frak
M$. Then $\frak F(\frak M) \neq \frak M$. Contradiction. Hence 1)
implies 2).

Conversely, let $\frak F(L) = L$ and we assume that $\frak F(\frak
M) \neq \frak M$. Let $\frak N$ be a maximal proper subgroup of
$\frak M$. $\frak M$ is a locally nilpotent group [2], then by
Lemma 4 the commutator subgroup $\frak M^{\prime}$ lies in $\frak
F(\frak M)$. Then $\frak N$ is a normal subgroup of $\frak M$ and
$\frak M/\frak N$ is a cyclic group of prime order $p$. By Lemma 1
$\frak N1 = H$ is a normal subloop of $L$ and $\frak N \subseteq
H^{\ast}$. Then from $\frak M(L/H) \cong \frak M/H^{\ast}$ it
follows that $L/H$ is a cyclic group of order $p$. Hence $H$ is a
maximal subloop of $L$. Then $\frak F(L) \neq L$. Contradiction.
Hence 2) implies 1). This completes the proof of Lemma 6.
\smallskip\\

Let $M(H)$ denote the subgroup of multiplication group of CML $L$,
gene\-ra\-ted by $\{L(x) \vert \forall x \in H\}$, where $H$ is a
subset of $L$.
\smallskip\\

\textbf{Lemma 7 [8].} \textit{Let $L$ be a CML with multiplication
group $\frak M(L)$ and inner mapping group $\frak I(L)$. Then
$\frak M(L)^{\prime} = <\frak I(L), M(L^{\prime})> =
(L^{\prime})^{\star} = \overline{\frak I(L)}$, where
$(L^{\prime})^{\star} = \{\alpha \in \frak M(L) \vert \alpha
x\cdot L^{\prime} = xL^{\prime} \quad \forall x \in L\}$,
$\overline{\frak I(L)}$ is the normal subgroup of $\frak M(L)$,
generated by $\frak I(L)$.}
\smallskip\\

\textbf{Proposition 2.} \textit{For a  CML $L$ with multiplication
group $\frak M$ the following statements are equiva\-lent: 1)
$\frak F(L) = L$ and $L$ satisfies the identity $x^3 = 1$; 2)  $L
= L^{\prime}$; 3) $\frak M = \frak M^{\prime}$; 4) $\frak F(L) =
L$ and $Z(L) = \{1\}$; 5) $\frak F(\frak M) = \frak M$ and
$Z(\frak M) = \{e\}$.}
\smallskip\\

\textbf{Proof}. 1) $\Rightarrow$ 2). As $\frak F(L) = L$ then by
Lemma 5 $\frak F(L/L^{\prime}) = L/L^{\prime}$. In [9] it is
proved that for an abelian group $G$ $\frak F(G) = G$ if and only
if $G$ is a divisible group. The abelian group $L/L^{\prime}$
satisfies the identity $x^3 = 1$. $L/L^{\prime}$ is a divisible
group, then $L/L^{\prime}$ is an unitary group. Hence $L^{\prime}
= L$, i.e. 1) implies 2). Conversely, let $L^{\prime} = L$. By [4]
the associator subloop $L^{\prime}$ satisfies the identity $x^3 =
1$ and from relations $L = L^{\prime} \subseteq \frak F(L)
\subseteq L$ it follows that $\frak F(L) = L$. Hence 2) implies
1).

1) $\Rightarrow$ 3). By Lemma 6 $\frak F(L) = L$ implies $\frak
F(\frak M) = \frak M$. Like in the previous case from here it
follows that $\frak M/\frak M^{\prime}$ is a divisible abelian
group. By definition the group $\frak M$ is generated by
translations $L(x), x \in L$. Then from identity $x^3 = 1$ for $L$
and di-associativity of $L$ it follows that the divisible abelian
group $\frak M/\frak M^{\prime}$ satisfies the identity $x^3 = 1$.
Then $\frak M/\frak M^{\prime}$ is an unitary group. Hence $\frak
M^{\prime} = \frak M$, i.e. 1) $\Rightarrow$ 3). Conversely, let
$\frak M^{\prime} = \frak M$. By Lemma 4 $\frak M^{\prime}
\subseteq \frak F(\frak M)$. Then from relations $\frak M = \frak
M^{\prime} \subseteq \frak F(\frak M) \subseteq \frak M$ it
follows that $\frak F(\frak M) = \frak M$. By Lemma 1 $\frak
M(L/L^{\prime}) \cong \frak M/(L^{\prime})^{\ast}$. $\frak
M(L/L^{\prime})$ is an abelian group. Then $\frak M^{\prime}
\subseteq (L^{\prime})^{\ast}$ and from the relation $\frak
M^{\prime} = \frak M$ it follows that $\frak M(L/L^{\prime})$ is
an unitary group. Hence $L^{\prime} = L$. Consequently, 3) implies
2).

2) $\Rightarrow$ 4). We consider the homomorphism $\alpha: L
\rightarrow L/Z(L)$. The elements of quotient loop have the form
$aZ(L)$, $a \in L$. From $L = L^{\prime}$ it follows that the
element $a$ is a product of associators $(u,v,w)$, $u, v, w \in
L$. From equalities $(u,v,w)Z(L) = (uZ(L),v,w) = (u,v,w)$,
$ab\cdot Z(L) = a\cdot bZ(L) = aZ(L)\cdot b$ it follows that if
$aZ(L) = bZ(L)$ then $a = b$. But this means that $\alpha$ is
isomorphism. Then $Z(L) = \{1\}$. Consequently, 2) implies 4).
Conversely, let $Z(L) = \{1\}$. Then from Lemma 2 it follows that
CML $L$ satisfies   the identity $x^3 = 1$. Hence 4) implies 1).
Further, the equivalence of statements 4), 5) follows from Lemma 6
and Proposition 1. This completes the proof of Proposition 2.
\smallskip\\

\textbf{Theorem 1.} \textit{For a CML $L$ with multiplication
group $\frak M$ the following statements are equivalent: 1) $\frak
F(L) = L$; 2) $L$ is a direct product $L = L^{\prime} \times L^3$,
where $L^{\prime}$ is the associator subloop of $L$ and $L^3 =
\{x^3 \vert x \in L\}$ is a divisible abelian group; 3) $\frak
F(\frak M) = \frak M$; 4) $\frak M$ is a direct product $\frak M =
\frak M^{\prime} \times \frak D$, where $\frak M^{\prime}$ is the
commutator subgroup of $\frak M$ and $\frak D$ is a divisible
abelian group. In such cases $\frak F(L^{\prime}) = L^{\prime}$,
$Z(L^{\prime}) = \{1\}$, $Z(L) = L^3$,  $L^3 \cong \frak D$,
$Z(\frak M) = \frak D$, $\frak F(\frak M^{\prime}) = \frak
M^{\prime} =  M(L^{\prime}) = (L^{\prime})^{\star} =
\overline{\frak I(L)}$, where $(L^{\prime})^{\star} = \{\alpha \in
\frak M(L) \vert \alpha x\cdot L^{\prime} = xL^{\prime} \quad
\forall x \in L\}$, $\overline{\frak I(L)}$ is the normal subgroup
of $\frak M(L)$, generated by inner mapping group $\frak I(L)$,
$Z(\frak M^{\prime}) = \{e\}$.}
\smallskip\\

\textbf{Proof.} 1) $\Rightarrow$ 2). Using the diassociativity of
CML it is easy to prove that $L^3$ is a subloop of $L$. By Lemma 2
$L^3 \subseteq Z(L)$. Then $L^3$ is a normal associative subloop
of $L$. The quotient loop $L/L^3$ satisfies the identity $x^3 =
1$. By Lemma 5 from $\frak F(L) = L$ it follows that $\frak
F(L/L^3) = L/L^3$. Then by Proposition 2 $L/L^3 =
(L/L^3)^{\prime}$. But $(L/L^3)^{\prime} = L^{\prime}L^3/L^3$.
Then from $L/L^3 = L^{\prime}L^3/L^3$ it follows that $L =
L^{\prime}L^3$. Hence $L/L^3 = L^{\prime}L^3/L^3 \cong
L^{\prime}/(L^{\prime} \cap L^3)$. We have $\frak
F(L^{\prime}/(L^{\prime} \cap L^3)) =
L^{\prime}/(L^{\prime}/(L^{\prime} \cap L^3))$ and $L^{\prime}
\cap L^3 \subseteq Z(L)$. Then by analogy with the proof of
implication 1) $\Rightarrow$ 4) of Proposition 2 it is easy to
prove that $L^{\prime} \cap L^3 = \{1\}$. But $L = L^{\prime}L^3$.
Then $L = L^{\prime} \times L^3$. Further, by Lemma 5 we get that
$\frak F(L^{\prime}) \cong \frak F(L/L^3) = \frak F(L)/L^3 = L/L^3
\cong L^{\prime}$, $\frak F(L^{\prime}) = L^{\prime}$ and by
Proposition 2 $Z(L^{\prime}) = \{1\}$. Analogously, $\frak F(L^3)
= L^3$. The subloop $L^3$ is associative. Then from $\frak F(L^3)
= L^3$ it follows that $L^3$ is a divisible abelian group [9].
Consequently,  1) implies 2) and $\frak F(L^{\prime}) =
L^{\prime}$, $Z(L^{\prime}) = \{1\}$. Further, from $L = L{\prime}
\times L^3$, $Z(L^{\prime}) = \{1\}$ it follows that $Z(L) = L^3$.

Conversely, let $L = L^{\prime} \times L^3$ and let $L^3 \subseteq
Z(L)$ be a divisible group. Then $\frak F(L^3) = L^3$ and
$L^{\prime} = (L^{\prime} \times L^3)^{\prime} =
(L^{\prime})^{\prime}$. By Proposition 2 $\frak F(L^{\prime}) =
L^{\prime}$. Hence $L = \frak F(L^{\prime}) \times \frak F(L^3)$.
$\frak F(L^{\prime})$ and $\frak F(L^3)$ do not have a maximal
proper subloops. From here it is easy to see that  $L$ doesn't
have a maximal proper subloops, either. Then $\frak F(L) = L$.
Hence 2) implies 1) and, consequently, the statements 1), 2) are
equivalent.

The equivalence of statements 1), 3) follows from Lemma 6.

2) $\Rightarrow$ 4). Let $L = L^{\prime} \times L^3$. From here it
follows that any element $a \in L$ has the form $a = ud$, where $u
\in L^{\prime}$, $d \in L^3$. As by Lemma 2 $L^3 \subseteq Z(L)$,
then $L(a) = L(u)L(d)$, therefore, $\frak M =
M(L^{\prime})M(L^3)$. Any element $\alpha \in M(L^3)$ has the form
$\alpha = L(v)$, where $v \in L^3$. Let $\alpha \in M(L^{\prime})
\cap M(L^3)$. Then $\alpha 1 \in L^{\prime} \cap L^3 = \{1\}$,
$\alpha 1 = 1$, $L(v)1 = 1$, $v = 1$, $L(v) = e$, $M(L^{\prime})
\cap M(L^3) = \{e\}$. Further, by Proposition 1 $Z(L) = L^3$
implies $Z(\frak M) = M(L^3) \cong Z(L) \cong \frak D$. $M(L^3)$
is a normal subgroup of $\frak M$. Then from $\frak M =
M(L^{\prime})M(L^3)$ it follows that  $M(L^{\prime})$ is also
normal in $\frak M$. Hence $\frak M = M(L^{\prime}) \times \frak
D$. The quotient loop $\frak M/M(L^{\prime})$ is abelian. Then
$\frak M^{\prime} \subseteq M(L^{\prime})$. By Lemma 7
$M(L^{\prime}) \subseteq \frak M^{\prime}$. Hence $M(L^{\prime}) =
\frak M^{\prime}$. Consequently, $\frak M = \frak M^{\prime}
\times \frak D$, i.e. 2) implies 4). Conversely, if $\frak M =
\frak M^{\prime} \times \frak D$ then $\frak M = M(L^{\prime})
\times M(L^3)$, $\frak M1 = M(L^{\prime})1 \times M(L^3)1$. Hence,
4) implies 2).

Finally, the equality $\frak F(\frak M^{\prime}) = \frak
M^{\prime}$ follows, by Lemma 5, from relations $\frak F(\frak M)
= \frak M$, $\frak M/M(L^3) \cong \frak M^{\prime}$, the
equalities $\frak M^{\prime} = (L^{\prime})^{\star} =
\overline{\frak I(L)}$ follow from Lemma 7 and the equality
$Z(\frak M^{\prime}) = \{e\}$ follows from equalities $\frak M =
\frak M^{\prime} \times M(L^3)$, $Z(\frak M) = M(L^3)$. This
completes the proof of Theorem 1.

\section{Normalizer condition}

Let $M$ be a subset, $H$ be a subgroup of group $G$. The subgroup
$N_H(M) \break = \{h \vert h \in H, h^{-1}Mh = M\}$ is called
\textit{normalizer} of set $M$ in subgroup $H$ [7]. Now,
constructively, we define  the notion of normalizer for the
subloops of CML. Let $H, K,$ where $H \subseteq K$ are subloops of
CML $L$. We define inductively the series of sets $\{P_{\alpha}\}$
and $\{D_{\alpha}\}$ as follows:

i) $P_1 = \{x \in K \vert (H,H,x) \subseteq H\}$ and $D_1 = \{x
\in K \vert (H,x,P_1) \subseteq H\}$;

ii) for any ordinal $\alpha$, $P_{\alpha+1} = \{x \in K \vert
(H,D_{\alpha},x) \subseteq H\}$ and $D_{\alpha+1} = \{x \in K
\vert (H,x,P_{\alpha+1}) \subseteq H\}$;

iii) if $\alpha$ is a limit ordinal, $P_{\alpha} = \bigcap_{\beta
< \alpha}P_{\beta}$ and $D_{\alpha} = \bigcup_{\beta <
\alpha}D_{\beta}$. \\ Further, we will also denote the conditions
of item ii) by  $(H,D_{\alpha},\overline P_{\alpha+1})$ and \break
$(H,\overline D_{\alpha+1},P_{\alpha+1})$ respectively. $H$ is a
subloop of CML $L$, then from \break $(H,H,H) \subseteq H$,
$(H,H,\overline P)$ it follows that $H \subseteq P_1$, from $H
\subseteq P_1$, $(H,\overline D_1,P_1)$, $(H,H,\overline P_1)$ it
follows that $H \subseteq D_1$, from $(H,H,\overline P_1)$,
$(H,D_1,\overline P_2)$, $H \subseteq D_1$ it follows that $P_1
\supseteq P_2$, from $(H,\overline D_1,P_2)$, $(H,\overline
D_2,P_2)$, $P_1 \supseteq P_2$ it follows that $D_1 \subseteq
D_2$. Further, let $\alpha$ be a non-limit ordinal and we suppose
by inductive hypothesis that $D_{\alpha} \supseteq D_{\alpha + 1}$
and $P_{\alpha} \subseteq P_{\alpha + 1}$. Then from $D_{\alpha}
\supseteq D_{\alpha + 1}$, $(H,D_{\alpha},\overline P_{\alpha +
1})$, $(H,D_{\alpha + 1},\overline P_{\alpha + 2})$ it follows
that $P_{\alpha + 1} \subseteq P_{\alpha + 2}$ and from $P_{\alpha
+ 1} \subseteq P_{\alpha + 2}$, $(H,\overline D_{\alpha +
2},P_{\alpha\alpha + 1})$, $(H,\overline D_{\alpha +
2},P_{\alpha\alpha + 2})$ it follows that $D_{\alpha + 1}
\supseteq D_{\alpha + 2}$. Hence, if consider also item iii), we
get a series of subsets
$$P_1 \supseteq P_2 \supseteq \ldots \supseteq P_{\alpha}
\supseteq \ldots$$
$$D_1 \subseteq D_2 \subseteq \ldots \subseteq D_{\alpha}
\subseteq \ldots. \eqno{(5)}$$

The construction process  of subsets $P_{\alpha}$, $D_{\alpha}$
from (5) shall end with an ordinal number, whose cardinality
doesn't exceed the cardinality of CML $K$ itself. We suppose that
$P_{\alpha + 1} = P_{\alpha + 2} = \dots.$ From $(H,\overline
D_{\alpha + 1},P_{\alpha + 1})$, $(H,\overline D_{\alpha +
2},P_{\alpha + 2})$ it follows that $D_{\alpha + 1} = D_{\alpha +
2}$. Then from $(H,D_{\alpha + 1},\overline P_{\alpha + 2})$,
$(H,\overline D_{\alpha + 2},P_{\alpha + 2})$ it follows that
$(H,\overline D_{\alpha + 2}, \overline P_{\alpha + 2})$. We
remind that the insc\-rip\-tions $\overline D_{\alpha + 2}$,
$\overline P_{\alpha + 2 }$ denote the biggest subsets $D_{\alpha
+ 2}$ and $P_{\alpha + 2}$ that especi\-ally the relation
$(H,D_{\alpha + 2},P_{\alpha + 2}) \subseteq H$ holds true. $H$ is
a subloop of CML $L$, then from (3) it follows that $D_{\alpha +
2} = P_{\alpha + 2}$ and using (3), (4) it is easy to prove that
$D_{\alpha + 2}$ is a subloop of CML $L$. Hence $D_{\alpha + 2}$
is the biggest (and the single) subloop of CML $K$ where  by (2)
$H$ is a normal subloop. By analogy with group theory the subloop
$D_{\alpha}$ will be called \textit{normalizer} of subloop $H$ in
subloop $K$ of CML $L$ and will be denoted by $N_K(H)$. If the
subgroup, where the normalizer is taken from, is not  indicated,
it means that it is taken from the entire  CML $L$. Consequently,
from construction of normalizer it follows.
\smallskip\\

\textbf{Proposition 3.} \textit{Let $H$, $L$, $K$, where $H
\subseteq L \subseteq K$, be  subloops of CML $L$ and let $H$ be a
normal subloop of $L$. Then $L \subseteq N_K(H)$.}
\smallskip\\

The group theory contains good studies of the group that satisfies
the normalizer condition (see, for example, [7]). These are such
groups, where every  proper subgroup differs from its normalizer.
A similar notion can be introduced for CML. We will say that a
\textit{CML satisfies the normalizer condition} or, in short, is a
$N$-loop, if every proper subloop differs from its normalizer.

\textit{The CML $L$ then and only then will be a $N$-loop, when an
ascending subnormal system $\{H_{\alpha}\}$  passes through each
subloop $H$ of CML $L$.}

Really, we denote $H_0 = 1, H_1 = H$ (respect. $\frak N_0 = e,
\frak N_1 = \frak N$). Further, for non-limit $\alpha$ we take as
$H_{\alpha}$  the normalizer of subloop $H_{\alpha-1}$, and for
limit $\alpha$ $H_{\alpha}$ will be  the union of all $H_{\beta}$
under $\beta < \alpha$. This ascending subnormal system,
obviously, reaches CML $L$  itself. Conversely, if all subloops of
CML $L$ are contained in some ascending subnormal system, then all
proper subloop will be normal in some bigger subloop, and,
consequently, by Proposition 3, will differ from its normalizer.

Using this result it is easy to prove that \textit{all subloops
and all quotient loops of $N$-loop will be $N$-loops themselves.}
Really, let $A$ be a subloop of $N$-loop $L$, and let $B$ be a
subloop of $L$ such that $B \subseteq A$. By the aforementioned,
an ascending subnormal system $\{B_{\alpha}\}$ passes through $B$.
Then $\{B_{\alpha} \cap A\}$ after removing the repetitions will
be an ascending subnormal system of $A$, passing through  $B$.
Hence $A$ will be a $N$-loop. The second statement is proved by
analogy.
\smallskip\\

\textbf{Theorem 2.}  \textit{If a CML $L$ with the Frattini
subloop $\frak F(L)$ satisfies the inequality $\frak F(L) \neq L$
then it satisfy the normalizer condition}.
\smallskip\\

\textbf{Proof.} As $\frak F(L) \neq L$ then  the CML $L$ has a
maximal proper subloops. Let $H$ be an arbitrary proper subloop of
CML $L$. If $H$ is a maximal subloop of $L$ then by [4] $H$ is
normal in $L$. Hence $H \neq N_L(H) = L$. Let now the subloop $H$
be a non-maximal subloop. By Zorn's Lemma let $M$ be a maximal
subloop of $L$ regarding the property $H \subseteq M$ and let $a
\notin L \backslash M$. We suppose that $a^3 = 1$. Let $K = <H,
a>$. $M$ is a maximal proper subloop of $L$, then by [4] the
subloop $M$ is normal in $L$.  Let $\varphi$ be a restriction in
$K$ of homomorphism $L \rightarrow L/M$. Obviously,  $\text{Ker}
\varphi = M \cap K$. As $a^3 = 1$ then $M \cap <a> = 1$. Hence $K
\backslash <a> = H$ and then $M \cap K = H$. Consequently, $H$ is
a normal subloop of $K$, and as $H \neq K$ then by Proposition 1
$H \neq Z_L(H)$, as required.

Let now $a^3 \neq 1$. By Lemma 1 $a^3 \in Z(L)$, hence $<a^3>$ is
a normal subloop of $L$. Let $a^3 \in H$. We denote $L/<a^3> =
\overline L$. From $a^3 \in M$, $a \notin M$ it follows that
$\overline M$ is a maximal proper subloop of $\overline L$. Hence
$\frak F(\overline L) \neq \overline L$. Further, $\overline a^3 =
\overline 1$, then the previous cases $\overline H \neq
N(\overline H)$. As $a^3 \in H$ and $a^3 \in N(H)$ then the
inverse images of $\overline H$ and $N(\overline H)$ will be $H$
and $N(H)$ respectively. Hence from $\overline H \neq N(\overline
H)$ it follows that $H \neq N(H)$, as required.

If $a^3 \notin H$, then $H \neq H<a^3>$. By (3) and Lemma 1 we get
$(H,H<a^3>,H<a^3>) = (H,H,H) \subseteq H$. This means by
Proposition 1 that $H<a^3> \subseteq N(H)$. Hence $H \neq N(H)$.
This completes the proof of  Theorem 2.
\smallskip\\

Any subloop of a $ZA$-loop is a $ZA$-loop. From Lemma 3 it follows
that a non-associ\-at\-ive commutative Moufang $ZA$-loop has a
non-trivial associati\-ve quotient loop. Hence it differs from its
associator subloop. Hence \textit{any commutative Moufang
$ZA$-loop is a $SD$-loop}.
\smallskip\\

\textbf{Corollary 1.} \textit{For a CML $L$ let $Z_2(L) \neq
Z_1(L)$. In particular, let $L$ be a $ZA$-loop or a SD-loop. Then
the CML $L$ satisfies the normalizer condition.}
\smallskip\\

\textbf{Proof.} We suppose that $\frak F(L) = L$. Then by Theorem
1 $L = L^{\prime} \times Z(L)$, $\frak F(L^{\prime}) =
L^{\prime}$, $Z(L^{\prime}) = \{1\}$. From here  it follows that
$Z_2(L) = Z_1(L)$. Contradiction. Hence $\frak F(L) \neq L$ and by
Theorem 2 the CML $L$ satisfies the normalizer condition, as
required.
\smallskip\\

Now we will enforce the Corollary 1 for centrally nilpotent CML
$L$ with multip\-lication group $\frak M$. In [4] it is proved
that CML $L$ is centrally nilpotent of class $n$ if and only if
the group $\frak M$ is nilpotent of class $2n - 1$. Then Corollary
1 for $\frak M$ follows from the known result about subnormal
subgroups of nilpotent group (see, for example, [7]).
\smallskip\\

\textbf{Proposition 4.} \textit{If $L$ is a centrally nilpotent
CML of class $n$, then for its any subloop $H$ (respect. subgroup
$\frak N$ of group $\frak M$) the series of consecutive
normalizers reaches $L$ (respect. $\frak M$) not later that after
$n$ (respect. $2n - 1$) steps.}
\smallskip\\

\textbf{Proof.} Let (1) be the upper central series of CML $L$. We
denote $H_0 = H$, $H_{i+1} = N_L(H_i)$. It is sufficient to check
that $Z_i(L) \subseteq H_i$. For $i = 0$, this is obvious. We
suppose that $Z_i(L) \subseteq H_i$. From the relation
$Z_{i+1}(L)/Z_i(L) = Z(L/Z_i(L)$ it follows that $(Z_{i+1}(L),L,L)
\subseteq Z_i(L)$. In particular, $(Z_i(L), \break
Z_{i+1}(L),Z_{i+1}(L)) \subseteq Z_i(L)$. As $Z_i(L) \subseteq
H_i$, then $(H_i,Z_{i+1}(L), Z_{i+1}(L)) \subseteq H_i$. But this
is to note down that $Z_{i+1}(L)$ normalizes $H_i$. Hence
$Z_{i+1}(L) \subseteq H_{i+1}$. This completes the proof of
Proposition 4.
\smallskip\\

\textbf{Remark}. Theorem 1 (see, also, Theorem 3) reveals a strong
analogy between Frattini subloops  of CML and Frattini subgroups
of multiplication groups of CML. However for multiplication group
of CML the statement, analogous to Theorem 2, is false. In [4]
there is an example of CML $G$ of exponent 3, such that
$G^{\prime} \neq G$ and $Z(G) = 1$. By Proposition 2 $\frak F(G)
\neq G$. Then by Proposition 1 $Z(\frak M) = e$ and by Lemma 6
$\frak F(M) \neq \frak M$, where $\frak M(G)$ denote the
multiplication group of $G$. In [5] J. D. H. Smith showed that no
group with trivial centre and satisfying the normalizer condition
can be the multiplication group of quasigroup. Hence the
multiplication group $\frak M$ satisfies the inequality $\frak
F(\frak M) \neq \frak M$ but it doesn't satisfy the normalizer
condition.

\section{Divisible subgroups of multiplication group}

We remind ([7]  (respect. [2])) that the group (respect. CML) $G$
is called \textit{divisible} or \textit{complete} (by terminology
of [3] \textit{radically complete}) if the equality $x^n = a$ has
at least one solution in $G$, for any number $n > 0$ and any
element $a \in G$.
\smallskip\\

\textbf{Theorem 3.}  \textit{Any divisible subgroup $\frak N$ of a
multiplication group $\frak M$ of a CML $L$ is an abelian group
and serves as direct factor for $\frak M$, i.e. $\frak M = \frak N
\times \frak C$ for a certain subgroup $\frak C$ of $\frak M$.}
\smallskip\\

\textbf{Proof.} If $\frak F(\frak M) = \frak M$ then the statement
follows from Theorem 1. Hence let $\frak F(\frak M) \neq \frak M$.
Then $\frak M$ has a maximal proper subgroups. The group $\frak M$
is locally nilpotent [2], then the maximal proper subgroups of
$\frak M$ are normal in $\frak M$ [7]. Let $\frak H$ be a maximal
proper subgroup of $\frak M$ such that $\varrho \notin \frak H$
for some $\varrho \in \frak N$. We will consider two cases:
$\varrho$ has a finite order and $\varrho$ has an infinite order.

Let the element $\varrho$ have a finite order $n$. Then the
element $\alpha = \varrho^{n/p}$, where $p$ is a prime divisor of
$n$, has the order $p$. The subgroup $\frak N$ is divisible. Then
there exists a series  $\alpha = \alpha_1, \alpha_2, \ldots,
\alpha_k, \ldots$ of elements in $\frak N$, such that $\alpha_1^p
= e$, $\alpha_{k +1}^p = \alpha_k$, where $e$ is the unit of
$\frak M$. From here it follows that $\alpha_k^{p^k} = e$, $k = 1,
2, \ldots$.

We denote by $\frak C$ the subgroup of $\frak N$ generate by
$\alpha_1, \alpha_2, \ldots, \alpha_k, \ldots$. It is easy to
prove that any element $\alpha \in \frak C$ is a power of some
generator $\alpha_k$, i.e. $\alpha = \alpha_k^n$, and the cyclic
groups $<\alpha_k>$ form a series
$$e \subset <\alpha_1> \subset <\alpha_2> \subset \ldots \subset
<\alpha_k> \subset \ldots.$$ We prove that $\frak C \cap \frak H =
e$. $<\alpha_1>$ is a cyclic group of order $p$ and $\alpha_1
\notin \frak H$. Then $<\alpha_1> \cap \frak H = e$. We suppose
that $<\alpha_k> \cap \frak H = e$. We have $<\alpha_{k+1}> =
\{\alpha_{k+1}, \alpha_{k+1}^p, \ldots, \alpha_{k+1}^{p^{k+1}-1}\}
\cup <\alpha_k>$. We suppose that $\alpha_{k+1}^n \in \frak H$ ($n
= 1, 2, \ldots, p^{k+1} - 1$). Then  $(\alpha_{k+1}^n)^p \in \frak
H$. But $(\alpha_{k+1}^n)^p = (\alpha_{k+1}^p)^n = \alpha_k^n$.
Hence $\alpha_k^n \in \frak H$. But this contradicts the
supposition $<\alpha_k> \cap \frak H = e$. Hence $<\alpha_{k+1}>
\cap \frak H = e$ and, consequently, $\frak C \cap \frak H = e$.
Let $\frak M^{\prime}$ denote the commutator subgroup of group
$\frak M$. By Lemma 4  $\frak M^{\prime} \subseteq \frak H$. Then
$\frak M^{\prime} \cap \frak C = e$,  $\frak C^{\prime} = e$,
hence $\frak C$ is an abelian group. More concretely, $\frak C$ is
isomorphic to a quasicyclic $p$-group. Further, the subgroup
$\frak H$ as maximal in $\frak M$ is normal in $\frak M$. Then
from $\frak M = \frak H\frak C$, $\frak H \cap \frak C = e$ it
follows that $\frak C$ is normal in $\frak M$. Hence $\frak M =
\frak H \times \frak C$.

Let now $\varrho \in \frak N$ be an element of infinite order. If
$Z(\frak M)$ denotes the centre of $\frak M$ then $\frak M/Z(\frak
M)$ is a locally finite $3$-group [4]. Hence $\varrho^n \in
Z(\frak M)$ for some $n$. Let $\frak H$ be a maximal subloop of
$\frak M$ such that $\varrho^n \notin \frak H$. $\frak N$ is a
divisible group. Then there exists a series $\varrho^n = \alpha_1,
\alpha_2, \ldots, \alpha_k, \ldots$ of elements in $\frak N$, such
that $\alpha_{k+1}^{k+1} = \alpha_k$, $k = 1, 2, \ldots$. We
denote by $\frak Q$ the subgroup of $\frak N$, generated by
$\varrho^n = \alpha_1, \alpha_2, \ldots, \alpha_k, \ldots$. As
$\alpha_1 \in Z(\frak M)$ then it is easy to see that $\frak Q
\subseteq Z(\frak M)$. Hence the subgroup $\frak Q$ is normal in
$\frak M$. The subgroup $\frak Q$ is without torsion. In [4] it is
proved that the commutator subgroup of multiplication group of any
CML is a locally finite $3$-group. Then $\frak Q \cap \frak
M^{\prime} = e$. From here it follows that $\frak Q^{\prime} = e$,
i.e. $\frak Q$ is an abelian group. More concretely, $\frak Q$ is
isomorphic to additive group of rationals.

Thus, in both  cases in group $\frak M$ there  exists an abelian
normal subgroup $\frak D \subseteq \frak N$, which is isomorphic
to quasicyclic $p$-group or additive group of rationals such that
$\frak M = \frak H \times \frak D$. We will use this procedure of
separating the  divisible subgroup from $\frak N$ as direct factor
for define the subgroups $\frak M_{\beta}$, $\frak A_{\beta}$ of
group $\frak M_{\beta - 1}$.

 Let $\frak M_0 = \frak M$, $\frak M_1 = \frak H$,
$\frak D_1 = \frak D$. For a non-limit ordinal $\beta$ inductively
we define $\frak M_{\beta - 1} = \frak M_{\beta} \times \frak
D_{\beta}$. We denote $\frak A_{\beta} = \frak D_1 \times \frak
D_2 \times \ldots \times \frak D_{\beta}$. As $\frak D_1, \frak
D_2, \ldots, \frak D_{\beta} \subseteq \frak N$ then $\frak
A_{\beta} \subseteq \frak N$. Further we consider the series of
subgroups $$\frak A_1 \subset \frak A_2 \subset \ldots \subset
\frak A_{\beta} \subset \ldots,$$
$$\frak M_1 \supset \frak M_2 \supset \ldots \supset \frak
M_{\beta} \supset \ldots, \beta < \alpha,$$ where $\frak M_{\beta
- 1} = \frak M_{\beta} \times \frak D_{\beta}$ if $\beta$ is a
non-limit ordinal and $\frak A_{\beta} = \cup_{\gamma <
\beta}\frak A_{\gamma}$, $\frak M_{\beta} = \cap_{\gamma <
\beta}\frak M_{\gamma}$ if $\beta$ is a limit ordinal.

It is clear that $\frak M_{\beta}$, $\frak A_{\beta}$ are normal
subgroups of $\frak M$. We prove that $\frak M = \frak M_{\beta}
\times \frak A_{\beta}$ for any $\beta$. If $\beta$ is a non-limit
ordinal, then by induction $\frak M = \frak M_1 \times \frak D_1 =
\frak M_1 \times \frak A_1 = \frak M_{\beta - 1} \times \frak
A_{\beta -1} = \frak M_{\beta} \times  \frak D_{\beta} \times
\frak A_{\beta -1} = \frak M_{\beta} \times \frak A_{\beta}$.
Hence $\frak M = \frak M_{\beta} \times \frak A_{\beta}$.

Let now $\beta$ be a limit ordinal and let $e \neq \lambda \in
\frak M_{\beta} \cap \frak A_{\beta}$. Then there exists a
non-ordinal $\delta < \beta$ that $\lambda \in \frak A_{\delta}$.
From $\lambda \in \frak M_{\beta} = \cap_{\gamma < \beta}\frak
M_{\gamma}$ it follows that $\lambda \in \frak M_{\gamma}$ for all
$\gamma < \beta$. But $\delta < \beta$. Then $\lambda \in \frak
M_{\delta} \cap \frak A_{\delta}$. Contradiction. Hence $\frak
M_{\beta} \cap \frak A_{\beta} = e$ and we may consider the direct
product $\frak M_{\beta} \times \frak A_{\beta}$.

Let $\lambda \in \frak M \backslash (\frak M_{\beta} \times \frak
A_{\beta})$. Then $\lambda \notin \frak M_{\beta}$, $\lambda
\notin \frak A_{\beta}$, i.e. $\lambda \notin \cap_{\gamma <
\beta}\frak M_{\gamma}$, $\lambda \notin \cup_{\gamma <
\beta}A_{\gamma}$. Hence $\lambda \notin \frak M_{\gamma}$ for all
$\gamma < \beta$ and from $\lambda \in \frak M$, $\frak M = \frak
M_{\gamma} \times \frak A_{\gamma}$ it follows that $\lambda \in
\cap_{\gamma < \beta}\frak A_{\gamma} = \frak A_{\beta}$. We wet a
contradiction. Hence $\frak M = \frak M_{\beta} \times \frak
A_{\beta}$ for all $\beta$.

The process of inductive construction $\frak A_{\alpha}$ will be
end on the first number $\gamma$, for which $\frak A_{\gamma} =
\frak N$. Consequently, $\frak M = \frak M_{\gamma} \times \frak
N$. This completes the proof of Theorem 3.
\smallskip\\

By Theorem 3 any multiplication group $\frak M$ of CML contains a
maximal divisible  associative subloop $\frak D$ and $\frak M =
\frak D \times \frak R$, where  obviously $\frak R$ is a
\textit{reduced CML}, meaning that it has no non-unitary divisible
subgroups. Consequently, we obtain
\smallskip\\

\textbf{Corollary 2.} \textit{Any multiplication group $\frak M$
of CML  $L$ is a direct product of the divisible abelian subgroup
$\frak D$ and the reduced subgroup $\frak R$. The sub\-group
$\frak D$ is unequivocally defined, the subgroup $\frak R$ is
defined exactly till the isomor\-phism.}
\smallskip\\

\textbf{Proof.} Let us prove the last statement. As $\frak D$ is
the maximal divisible subgroup  of the multiplication group $\frak
M$, it is entirely characteristic in $\frak M$, i.e. it is
invariant regarding the endomorphisms of the group $\frak M$. Let
$\frak M = \frak D' \times \frak R'$, where $\frak D'$ is a
divisible subgroup, and $\frak R'$ is a reduced subgroup of the
group $\frak M$. We denote by $\varphi, \psi$ the endomorphisms
$\varphi:\frak M \rightarrow \frak D', \psi: \frak M \rightarrow
\frak R'$. As $\frak D$ is an entirely characteristic subgroup,
then $\varphi \frak D$ and $\psi \frak D$ are subgroups of the
group $\frak M$. It follows from the inclusions $\varphi \frak D
\subseteq \frak D'$ and $\psi \frak D \subseteq \frak R'$ that
$\varphi \frak D \cap \psi \frak D = 1$. By Theorem 3 $\frak D$ is
a abelian group, therefore $\varphi \frak D, \psi \frak D$ are
normal in $\frak D$. Then $d = \varphi d\cdot \psi d$ ($d \in
\frak D$) gives $\frak D = \varphi \frak D \cdot \psi \frak D$, so
$\frak D = \varphi \frak D \times \psi \frak D$. Obviously,
$\varphi \frak D \subseteq \frak D \cap \frak D', \psi \frak D
\subseteq \frak D \cap \frak R'$, then $\varphi \frak D = \frak D
\cap \frak D', \psi \frak D = \frak D \cap \frak R'$. Hence $\frak
D = (\frak D \cap \frak D')\times (\frak D \cap \frak R')$. But
$\frak D \cap \frak R' = 1$ as a direct factor of the divisible
group, that is contained by the reduced group. Therefore, $\frak D
\cap \frak D' \subseteq \frak D, \frak D \subseteq \frak D'$, i.e.
$\frak D = \frak D'$. This completes the proof of Corollary 2.
\smallskip\\

\smallskip

Tiraspol State University, Moldova

The autor's home address:

Deleanu str 1, Apartment 60

Chishinev MD-2071, Moldova

e-mail: sandumn@yahoo.com
\end{document}